\documentclass[amstex,amstext,amsfonts,epsf,reqno,11pt] {amsart}
\usepackage{epsfig,latexsym,amsfonts,color,
amsmath,amssymb,dsfont}
\usepackage[english]{babel}
\usepackage[a4paper, total={5.7in, 8in}]{geometry}
\def\eps{{\epsilon}}

\def\NT{\vert\!\vert\!\vert}

\def\Wo{{\stackrel{\circ}{W}}{}}

\def\IntO{\int\limits_{\Omega}}
\def\IntQ{\int\limits_{Q_T}}
\def\IT0{\int\limits^T_0}

\def\dvg{{\rm div}}
\def\pz{p^*}
\def\yz{y^*}

\def\ez{e^*}

\def\res{{\mathbf r}}
\def\a{{\rm a}}

\def\Rd{{\mathbb R}^d}
\def\Md{{\mathbb M}^{d\times d}}
\def\measure{{\boldsymbol \mu}}
\def\Am{{A^{-1}}\!}
\newcommand\be{\begin{eqnarray*}}
\newcommand\ee{\end{eqnarray*}}
\newcommand\ben{\begin{eqnarray}}
\newcommand\een{\end{eqnarray}}
\numberwithin{equation}{section}
\newtheorem{Theorem}{Theorem}

\newtheorem{Remark}{Remark}

\begin{document}
%\large
%\setcounter{section}{-1}
\title
{A posteriori estimates for inverse evolutionary
 problems}
\author{S. REPIN }
\address{ St. Petersburg Department of V.A. Steklov Institute of Mathematics, Fontanka 27, 191023, St. Petersburg, Russia}
\email{repin@@pdmi.ras.ru}
\subjclass{Primary 65N30}
\keywords{parabolic equations, ill-posed problems, inverse problems,
deviations from exact solution, error identities
a posteriori estimates of the functional type}
\maketitle

%%%%%%%%%%%%%%%%%%%%%%%%%%%%%%%%%%%%%
\begin{quote}
\small
The paper is focused on a posteriori error analysis of
 evolutionary problems,
 which are ill-posed in the sense of Hadamard. They may have
 no solutions and respective approximations may be highly unstable.
 For these problems, the a posteriori error estimation
 concept applicable to  well-posed problems is not valid and requires a significant
  reformulation. 
 A posteriori estimates of the functional type provide the basis 
 for addressing this issue. They 
  generate computable quantities
 that estimate the distance between a function (approximate solution to the inverse
 problem) and exact solution of a "close" well-posed
 initial boundary value problem associated 
 with the same differential equation. The solution of this close problem exactly satisfies one part of the imposed
 conditions  and the others  are satisfied approximately  with the accuracy $\eps$. This idea is elaborated for
 several evolutionary problems generated by the reaction-convection-diffusion
 equation: backward parabolic problem,
 problem with non-local time conditions, Cauchy parabolic problem, and parameter(s)
 identification problem. The estimates derived not only 
 control the accuracy of approximate solutions, 
 but also suggest a unified approach to calculating them numerically.
 
\end{quote}
%%%%%%%%%%%%%%%%%%%%%%%%%%%%%%%%%%%%%%%%%%%%%%%%%%%%%%%%%%%%%
\section{Introduction}
\subsection{A class of inverse problems}
Foundations of the theory of inverse problems were laid quite 
a long time ago.
Currently, this theory covers a large number of mathematical models that 
are not well-posed in the sense
  of Hadamard (e.g., see \cite{KabLav,LattLion,TikhAr}). In this paper, 
  we are concerned with
  one class of these problems generated by
 evolutionary equations of the type
\ben
\label{bwmod1}
u_t+A(t)u=0.
\een
In \cite{LattLion}, the authors study them assuming that
 the operator $A$ is generated by the bilinear form $a(t;u,v)$
defined on a Hilbert space $V$,
which is $V$--elliptic for all $t$. 
Let $H$ be another Hilbert space
 such that $V\subset H$ and $V$ is dense in $H$.
Then for any $\xi\in H$ there exists a 
unique solution of (\ref{bwmod1}) such that $u(0)=\xi$,
moreover $u\in L_2(0,T;V)$ and continuously depends
 on $\xi$.  This  problem  (which is called "forward" or "direct")
 is well-posed according to Hadamard's definition.
  
Another (backward)  problem is to find $\xi$ such that
$u(T,\xi)=\chi(x)$, where $\chi(x)\in H$ 
is a certain given function.
In \cite{LattLion}, it is reformulated as an 
optimal control problem: find $u\in V$ satisfying (\ref{bwmod1}) that
minimizes  the cost functional
\ben
\label{bwmod2}
J(\xi):=\IntO|u(x,T;\xi)-\chi(x)|^2dx.
\een
For this class of problems, 
\ben
\label{bwmod3}
\inf\limits_{\xi\in H}J(\xi)=0.
\een
It is well known that backward problems of this type 
are ill-posed and have a solution provided that certain
additional conditions are imposed on the function $\chi$ (e.g., analyticity).
These solutions are highly unstable, which generates serious difficulties
in quantitative analysis. 
Hence even if (\ref{bwmod3}) is proved,
it may be very difficult to construct a reliable approximation of the inverse
problem.
Known approaches to the construction
of approximations to backward parabolic problems (e.g., the "quasi--reversibility
method" \cite{LattLion} and nonlocal regularisation method \cite{Show})
use regularised problems with small parameter $\eps$ at the higher order differential
operator or in the initial condition. Theoretically, they generate
a minimising sequence for the functional (\ref{bwmod2})
if $\eps\rightarrow 0$. However,
 getting sharp approximations of  regularised problems with guaranteed small
errors becomes a difficult task.
Similar difficulties and analogous questions arise when studying 
ill-posed problems of other types.
The question that arises is how to define the notion  "reliable approximation"
to an inverse problem.

\subsection{A posteriori analysis of inverse and ill-posed problems}
Direct and inverse problems are very different, so that the
concept of a posteriori
error control for them should be formulated differently.
Indeed, a well-posed problem  is uniquely solvable.
Let $u$ be the solution and $v$ be an approximation found by
some numerical procedure. The question
\begin{quote}
 {\em Which neighbourhood of the   solution $u$  does an approximation $v$
belong to?}
\end{quote}
 is quite logical. A posteriori estimates of the functional type
  give a constructive answer
 to this question for any $v$ in the respective energy space
  (see \cite{ReGruyter}, where the reader will find
 a systematic discussion of these estimates
 for a wide spectrum of problems).

An ill-posed problem may have
no solution so that the question posed above is formally incorrect.
It clearly needs  a  modification, which takes into account specific 
features of inverse
problems where the notion of "solution" is often 
replaced by a "quasi-solution". The latter  may be understood as a solution 
to some similar (or close) well-posed problem. 
Therefore, it is logical to
reformulate the error control concept accordingly.
Let $v_\eps$ be a given function in the admissible functional
class viewed as an approximation  
of an inverse (ill-posed) problem.
We wish to know
\begin{quote}
{\em Which neighbourhood of $v_\epsilon$ contains the exact solution 
of a nearby well-posed problem?}
 \end{quote} 
Below
we suggest a method able to answer this question for a class of inverse parabolic
problems. It defines
 an $\eps$-neighbourhood of $v_\eps$  that contains the solution
of a direct (well-posed)  problem (where $\eps$
is explicitly known). This solution exactly satisfies one part
of the conditions while the other(s) are satisfied approximately with the accuracy $\eps$.
Thus, we can point out a well-posed problem, whose solution
can be considered as an
$\eps$-solution to the problem in question.

\subsection{Illustration of the  main idea}
The approach discussed in this paper is
based on the theory of {\em functional a posteriori estimates} that has
been earlier developed for direct boundary value problems (see \cite{ReGruyter,Re2021}
and references cited therein).
The main idea   can be explained 
with the paradigm of a simple parabolic problem (\ref{dir1}) -- (\ref{dir3}). 

Let $\Omega$ be a bounded domain in $\Rd$,
$d\in \{1,2,3\}$ with Lipschitz boundary $\Gamma$. The direct evolutionary problem is to find
$u(x,t)$ in the space--time cylinder $Q_T:=\Omega\times(0,T)$, $T>0$,
such that
\ben
\label{dir1}
&u_t-\Delta u=f\qquad&\text{in}\;Q_T,\\
\label{dir2}
&u(x,t)=0\qquad&\text{on}\;S_T:=\Gamma\times(0,T),\\
\label{dir3}
&u(x,0)=\xi(x).
\een
 It is well known that under standard assumptions on
$f$ and $\xi$ (e.g., $\int^T_0 \|f\|^2_\Omega dt<\infty$ and $\xi\in L_2(\Omega)$) the problem (\ref{dir1})--(\ref{dir3}) is well-posed
and the function $u(x,t)$ exists and is unique.

%supplied with the norm 
%$\|\yz\|_{\dvg,Q_T}:=\left(\|\yz\|^2_{Q_T}+\|\dvg \yz\|^2_{Q_T}\right)^{1/2}$.
Let $v$ and $\yz$ be two functions considered as
 approximations of $u$ and $\pz:=\nabla u$,
respectively. Assume that
\be
v\in V_0:=\bigl\{v\in V:=H^1(Q_T)\,\mid\;
\;v(x,t)=0\;\text{on}\;S_T\bigr\}
\ee
and 
% TODO (B6): original Y^*(Q_T):=L_2(\Omega,\Md) looked dimensionally inconsistent
% with a space-time flux; left unchanged pending your confirmation (likely L_2(Q_T,\Md)).
\be
\yz\in Y^*_\dvg:=\bigl\{y^*\in Y^*(Q_T):=L_2(\Omega,\Md)\,\mid\,{\rm div} y^*\in
L_2(Q_T)\bigr\}.
\ee
In \cite{ReLincei,Re2021}, it was proven that the errors $e:=v-u$
and $\ez:=\yz-\pz$ meet the following estimate.
\begin{Theorem}
\label{th1}
For any $v\in V_0 $, 
$\yz\in Y^*_\dvg(Q_T)$, and
$\beta\in (0,1]$ it holds
% TODO (B2): argument order of \res differs from its definition below.
% Keep one order consistently.
\begin{multline}
\label{er1}
(1-\beta)\|\nabla e\|^2_{Q_T}+
\|e^*\|^2_{Q_T}+\|e(\cdot,T)\|^2_\Omega\\
\leq \|e(\cdot,0)\|^2_\Omega+
\|\nabla v-\yz\|^2_{Q_T}+{\beta}^{-1}C^2_\Omega\|\res(v,\yz)\|^2_{Q_T},
\end{multline}
where 
$$
\res(\yz,v):=\dvg\yz+f-v_t
$$ 
and $C_\Omega$ is a constant in the Friedrichs inequality
for the domain $\Omega$.
\end{Theorem}

Consider now the backward evolutionary problem
\ben
\label{sim1}
&u_t-\Delta u=f\qquad&\text{in}\;Q_T,\\
\label{sim2}
&u(x,t)=0\qquad&\text{on}\;S_T,\\
\label{sim3}
&u(x,T)=\chi(x)\in\, L_2(\Omega).
\een
Assume that by some method we have found the functions $\yz_\eps\in Y^*_\dvg$  
and $v_\eps\in V_{0}$ such that
$v_\eps(x,T)=\chi(x)$ for all  $x\in \Omega$ and
\ben
\label{Iopt}
\Upsilon_\beta(v_\eps,\yz_\eps)\,\leq\,\eps^2,
\een
where
\be
\Upsilon_\beta(v,\yz):=\|\nabla v-\yz\|^2_{Q_T}+
{\beta}^{-1}C^2_\Omega\|\res(v,\yz)\|^2_{Q_T}.
\ee
For $\xi(x)=v_\eps(x,0)$ the problem (\ref{dir1})--(\ref{dir3}) has a solution $u_\eps\in V_0$.
We can view $v_\eps$ and $\yz_\eps$ as
approximations of $u_\eps$ and $\pz_\eps:=\nabla u_\eps$, respectively. 
By (\ref{er1}) we estimate the distance between $(v_\eps,\yz_\eps)$ and $(u_\eps,\pz_\eps)$.
We use (\ref{Iopt}) and notice that 
$$e(\cdot,0)=v_\eps(\cdot,0)-u_\eps(\cdot,0)=0.
$$
Therefore, (\ref{er1}) reads as follows:
\ben
\label{er2}
(1-\beta)\|\nabla (v_\eps-u_\eps)\|^2_{Q_T}+
\|\yz_\eps-\nabla u_\eps\|^2_{Q_T}+\|u_\eps(\cdot,T)-\chi\|^2_\Omega\;
\leq\eps^2.
\een
In other words, $\eps$-neighbourhoods of $v_\eps$ and $\yz_\eps$ 
contain the exact solution
$u_\eps(x,t)$
of a direct evolutionary problem with the initial condition $v_\eps(x,0)$ and the corresponding
flux $\pz_\eps(x,t)$. Also, we can
guarantee that
\ben
\label{er3}
\|u_\eps(\cdot,T)-\chi\|_\Omega\leq \eps.
\een
% TODO (B5): original referred to (\ref{bwmod2})--(\ref{bwmod3}) (the optimal-control
% formulation). Likely should be (\ref{sim1})--(\ref{sim3}). Left unchanged.
Hence $u_\eps$ is an $\eps$-solution of the problem (\ref{bwmod2})--(\ref{bwmod3}).
 Moreover,  
 direct minimization of the functional $\Upsilon_\beta(v,\yz)$ with respect
 to its arguments suggests a method of getting successful
 approximations of the ill-posed problem (this question is briefly discussed later, see  Section \ref{num}).

\subsection{Outline of the paper}
The  paper is organised as follows.
In Section 2, we briefly discuss a class of direct initial boundary value
problems associated with  the reaction--convection--diffusion equation.
Theorems \ref{th2} and \ref{th3} of Section 3 establish basic error
relations for the direct evolutionary problems. 
Here we apply the derivation method used in (\cite{Re2021,Re2022})
and extend it to problems with mixed  boundary conditions.
The results are used in the next section where
we consider the backward parabolic problem, problem with non-local time conditions,
initial--boundary value problem with data given
on a part of the boundary,
and discuss parameter identification problems. 
 Using functional
relations derived in Section 3, we deduce a posteriori estimates for these generally
ill-posed problems. The corresponding results are presented in
Theorems \ref{th4}--\ref{th6}. Using them we 
define a direct evolutionary problem that can be viewed as an $\eps$-solution
 to the ill-posed problem under consideration.
 
%%%%%%%%%%%%%%%%%%%%%%%%%%%%%%%%%%%%%%%%%%%%%%%%

\section{A class of initial boundary value problems}

Henceforth,
we assume that $\Gamma$ is disjointly decomposed into a Dirichlet part $\Gamma_1$
and a Neumann part $\Gamma_2$.
Consider the problem
\ben
\label{rd1}
&&u_t-\dvg \pz+\a\cdot\nabla u+\rho^2  u=f\;{\rm in}\;Q_T,\\
\label{rd2}
&&\pz=A\nabla u,\\
\label{rd3}
&&u(x,t)=0\qquad{\rm on}\;S_{1T}:=\Gamma_1\times (0,T),\\
\label{rd4}
&& \pz\cdot n=F\qquad{\rm on}\;S_{2T}:=\Gamma_2\times (0,T),\\
\label{rd5}
&&u(x,0)=\xi.
\een
Here, $u_t$ is the time derivative of $u(x,t)$, $n$ is the unit outward normal to $\Gamma$, $\rho(x)$ and $\a(x)$ are the reaction and convection parameters, respectively,
which satisfy the following conditions:
\ben
\label{rd6}
&&{\a}\in L_\infty(\Omega,\Rd),\;
\dvg \a\in L_\infty(\Omega),
\\
&&
\label{rd7}
\rho\in L_\infty(\Omega),\quad 0\leq\rho_\ominus\leq
\rho \leq\,\rho_\oplus,\\
\label{rd8}
&& \varkappa^2:=  -\frac{1}{2}\dvg  {\a}+\rho^2
\geq 0.
\een
It is assumed that the boundary $\Gamma$ is Lipschitz so that the unit outward normal vector $n$
is uniquely defined at almost all points of $\Gamma$.
Also, we assume that the inflow part of the boundary where $(\a\cdot n)(x)<0$
is a subset of $\Gamma_1$, i.e.,
\ben
\label{rd9}
\delta:=\a\cdot n\;\geq 0\;\;\;\;\text{on}\;\;\Gamma_2.
\een
%\rc{Sign convention check: with $n$ outward, $(\a\cdot n)<0$ means 
%inflow into $\Omega$; the two sentences must be reconciled explicitly.\\
%*** No contradiction. If all points where $\a\cdot n<0$ are on $\Gamma_1$ then on $\Gamma_2$
%holds (\ref{rd9})}

The problem has two important special cases: the reaction-diffusion 
($\a = 0$)
 and convection-diffusion ($\rho=0$) problems.
 
Coefficients of the symmetric matrix $A$ are bounded,
do not depend on $t$, and satisfy the condition
\ben
c^2_1| \xi|^2\leq A \xi\cdot  \xi\leq c^2_2|\xi|^2\qquad\forall\xi\in\Rd,\;
c_1>0.
\label{rd10}
\een
Also, we assume that
\footnote{
The method and corresponding results presented in the paper 
 can be extended to the problems
where $A$, $a$, and $\rho$  are 
 bounded functions of time.
 Here we avoid these generalizations, which 
 are not essential in the context of the main result presented in Section 4.}
\ben
\label{rd11}
  f~\in ~L_2(Q_T),\qquad 
\xi\in \Wo^1_2(\Omega), \qquad F\in L_2(S_{2T}). 
\een

Throughout the paper, we  use standard notation
$L_p(\Omega)$
and $W^l_p(\Omega)$  (where $l,p\geq 1$)
for the Lebesgue and Sobolev spaces, respectively. A space
is marked above by $\circ$ if the respective functions
vanish on $S_{1T}$. 

By  $\|\cdot\|_\Omega$ we denote the 
 $L_2$ norm of a function
in $\Omega$ and $\|\cdot\|_{Q_T}$ denotes analogous
norm in the space--time cylinder.  
Equivalent seminorms generated by $A$ and $A^{-1}$
are denoted by  $\|\cdot\|_{A,\Omega}$ and  $\|\cdot\|_{A^{-1},\Omega}$,
respectively, i.e., 
$$\|w\|_{A,\Omega}:=\left(\int_\Omega A\nabla w\cdot\nabla w dx\right)^{1/2}\quad\text{and}\quad
 \|w\|_{A^{-1},\Omega}:=\left(\int_\Omega A^{-1}\nabla w\cdot\nabla w dx\right)^{1/2}. 
 $$
 The notation
 $\|\cdot\|_{A,Q_T}$ and  $\|\cdot\|_{A^{-1},Q_T}$ is  used for functions in $Q_T$.

In what follows, we use the spaces
\be
W^{1,0}_2(Q_T):=L_2(0,T;W^{1}_2  (\Omega))\;{\rm and}\;
\stackrel{\circ}{W}\!\!^{1,0}_2(Q_T):=L_2(0,T; \stackrel{\circ}{W}\!\!^1_2  (\Omega))
\ee
supplied with the norm 
\be
\|w\|_{1,Q_T}:=\left(\int\limits^T_0(\|\nabla w\|^2_\Omega
+\|w\|^2_\Omega)\, dt\right)^{1/2}
\ee
and the space
\ben
\label{Vnul}
\stackrel{\circ}{W}\!\!^{1,1}_2(Q_T):=\left\{w\in W^{1,1}_2(Q_T),\;w=0\;
{\rm on}\;S_{1T}\right\}. 
\een
For the latter space we also use the abridged notation $V_0(Q_T)$.

For the  functions in $\stackrel{\circ}{W}\!\!^{1,0}_2(Q_T,\Rd)$, the norm 
$\|\nabla w\|_{Q_T}$ is equivalent to  $\|w\|_{1,Q_T}$. 
Here, and later on the symbol := means "equals by definition".
%
%For functions in this space we introduce
%weighted norms ($\mu,\,\nu>0$)
% \ben
% \label{nnn}
%\NNN w\NNN_{\mu,\nu,Q_T}:=\Bigl(\mu\|\nabla w\|^2_{Q_T}+\nu\|w_t\|^2_{Q_T}\Bigr)^{1/2}.
%\een

Following \cite{La}, by $\Wo^{\Delta,1}_{2}(Q_T)$ we denote the subspace
of $\Wo^{1,0}_2(Q_T)$
that consists of the
functions such that
\be
\|w\|_{\Delta,1,Q_T}:=
\int\limits_{Q_T}(w^2+w^2_{t}+\mid\nabla w\mid^2+
(\Delta w)^2)\,dx\,dt<+\infty
\ee
and the space
$V_2(Q_T)$, which is the Banach space
of functions from $W^{1,0}_2(Q_T)$
having finite norm
\be
\NT w \NT^2:=
{\rm vrai}\!\max\limits_{t\in (0,T)}
\|w(\cdot,t)\|^2_{\Omega} +\|\nabla w\|^2_{2,Q_T}.
\ee
Spatial gradient and divergence are denoted by $\nabla$ and $\dvg$,
respectively. For a function $w\in H^1(\Omega)$
vanishing on $\Gamma_1$ it holds  the Friedrichs inequality. A constant
in this inequality is denoted by $C_\Omega$.
Then,
\ben
\label{FrA}
\|w\|_\Omega\;\leq\;C_{A,\Omega} \|\nabla w\|_{A,\Omega},
\een
where $C_{A,\Omega}\leq\,c^{-1}_1{C_\Omega}$.

Also, we use functional spaces  associated with vector valued functions
(fluxes). They are
$Y^*(Q_T):=L_2(Q_T,\Rd)$ and the space
\ben
\label{Ydiv}
 Y^*_\dvg(Q_T):=\bigl\{y^*\in Y^*(Q_T)\,\mid\,{\rm div} y^*\in
L_2(Q_T),\;y^*\cdot n\in L_2(S_{2T})\;\text{\rm on}\;S_{2T}\bigr\}.
\een
%supplied with the norm 
%$\|\yz\|_{\dvg,Q_T}:=\left(\|\yz\|^2_{Q_T}+\|\dvg \yz\|^2_{Q_T}\right)^{1/2}$.
%%%%%%%%%%%%%%%%%%%%%%%%%%%%%%%%%%%%%
Generalized solution of (\ref{rd1})--(\ref{rd5})
can be defined as the function
$$
u\in V^{1,0}(Q_T):=C([0,T];L_2(\Omega))\cap W^{1,0}_2(Q_T)
$$
that satisfies the
integral identity (see
\cite{LaSoUr,La})
\begin{multline}
\int\limits_{Q_T}\Bigl(A \nabla u\cdot \nabla  w+({\rm a}\cdot\nabla u)w+\rho^2 u w\Bigr)
\,dx\,dt - \int\limits_{Q_T}u w_{t}\,dx\,dt\\
+
\int\limits_{\Omega}\bigl(u(x,T) w(x,T)-u(x,0) w(x,0)\bigr)dx\\
=
\int\limits_{Q_T}f w \,dx\,dt+\int\limits_{S_{2T}}F w dsdt
\quad\forall w\in V_0(Q_T)
\label{rdint}
\end{multline}
for any test function $ w\in V_0(Q_T)$.
Under the conditions (\ref{rd6})--(\ref{rd11}), the solution is more regular and
belongs to $\Wo^{\Delta,1}_2(Q_T)$.
%\bc{$f\in C(0,T;L_2(\Omega))$, $a\in C(0,T;W^{1,\infty}(\Omega,\Rd)$,
%$b\in C(0,T;L_\infty(\Omega)$}

 %%%%%%%%%%%%%%%%%%%%%%%%%%%%%%%%%%%%%%%%%%%%%%%%%%%
\section{Error identities for the problem (\ref{rd1})--(\ref{rd5})}
\label{ApoIdn}

In this section, we deduce error relations that hold for the direct initial boundary value problem
(\ref{rd1})--(\ref{rd5}). Let $u\in V_0(Q_T)$ and let two functions
(considered as
  approximations of
$u(x,t)$ and $\pz(x,t):=A\nabla u$) satisfy the conditions, 
 $$v(x,t)\in  V_0(Q_T)\quad\text{and}\quad\yz(x,t)\in Y^*_\dvg(Q_T).
 $$
 In other words, minimal regularity conditions are imposed on the approximations, 
 and no special
properties of these functions (e.g., Galerkin orthogonality) are required.
 The functions
$$e:=v-u\quad\text{and}\quad   \ez:=\yz-\pz
$$
are the corresponding errors.  The  functions
\ben
\label{resf}
 && \res_f(v,\yz):=f-v_t+\dvg \yz-\a\cdot\nabla v-\rho^2  v, \\
 \label{resF}
 &&\res_F(\yz)=F-y^*\cdot n,\\
 \label{resA}
  &&\res_A(v,\yz):=A\nabla v-\yz. 
\een
contain
 known functions $v$ and $\yz$ and 
known data $f$, $F$, $A$, $\a$, and $\rho$. Therefore, 
they are directly computable and represent residuals
of the relations (\ref{rd1}),  (\ref{rd2}), and (\ref{rd4}). 
%%%%%%%%%%%%%%%%%%%%%%%%%%%%%%%%%%%%%%%%%%%%%%%%%%%%%%%%%%
\subsection{The main error identity}
The quantity $\measure(e,\ez)$ defined by the relation
\be
\measure^2(e,\ez):=
\|\nabla e\|^2_{A,Q_T}+
\|e^*\|^2_{\Am,Q_T}+\|e(\cdot,T)\|^2_\Omega
\ee
is a natural measure of the distance between the pair $(u,p^*)$ and the pair
$(v,\yz)$.  First, we establish identity (\ref{IDmain}), where $\measure(e,\ez)$ (plus some nonnegative
terms) form  the left hand side and the quantities $ \res_f(v,\yz)$ $ \res_F(v,\yz)$, and
$ \res_A(v,\yz)$ appear at the right  one. This identity plays a key role in further analysis.

%%%%%%%%%%%%%%
 \begin{Theorem}
 \label{th2}
 For any $v\in V_0(Q_T)$ and $y^*\in Y^*_\dvg(Q_T)$, it holds
 \begin{multline}
\label{IDmain}
\measure^2(e,\ez)+2\|\varkappa e\|^2_{Q_T}+\int\limits_{S_{2T}}(\a\cdot n)e^2dsdt\\
=\|v(\cdot,0)-\xi\|^2_\Omega+
\|\res_A(v,\yz)\|^2_{\Am,Q_T}\!\!-2\IntQ \!\!\res_f(v,\yz) e\; dxdt-
2\int\limits_{S_{2T}}\!\!\res_F(\yz)e dsdt.
\end{multline}
 \end{Theorem}
 \begin{proof}
We have
\begin{multline}
\label{adv}
\|\res_A(v,y^*)\|^2_{\Am,Q_T}=\|A\nabla e-\ez\|^2_{\Am,Q_T}\\=
\|\nabla e\|^2_{A,Q_T}
+\|\ez\|^2_{\Am,Q_T}-2\IntQ\ez\cdot\nabla e dxdt
\end{multline}
 and
 \begin{multline}
 \label{ba1}
 2\IntQ \!\!\res_f(v,\yz) e\; dxdt=
 2\IntQ(f-v_t+\dvg \yz-\a\cdot\nabla v-\rho^2 v)e\,dxdt\\
 =2\IntQ(\dvg \ez-e_t-\a\cdot\nabla e-\rho^2 e)\,e\,dxdt.
 \end{multline}
 By the relations
\ben
\label{ba2}
 \int\limits_{Q_T}e_t e\,dx\,dt=\frac12\IntO (e^2(x,T)-e^2(x,0))dx
\een
and
\ben
\label{ba3}
&&\IntO (\a\cdot\nabla e)edx=
\frac12\IntO \a\cdot\nabla (e^2)dx=
\frac12\int\limits_{\Gamma_{2}}(\a\cdot n)e^2ds-\frac12\IntO \dvg \a \,e^2 dx,
\een
{(the integral over $\Gamma_{1}$ vanishes since $e=0$ there)},
we rewrite (\ref{ba1}) as follows:
 \begin{multline}
 \label{ba4}
 2\IntQ \!\!\res_f(v,\yz) e\; dxdt=2\IntQ(e\dvg \ez-(\rho^2-\frac12\dvg \a)e^2)dxdt-
 \int\limits_{S_{2T}}(\a\cdot n)e^2dsdt\\
 +\IntO \bigl(e^2(x,0)-e^2(x,T)\bigr)dx.
 \end{multline}
 Now (\ref{adv}) and (\ref{ba4}) imply
 \begin{multline}
 \label{ba5}
\|\res_A(v,\yz)\|^2_{\Am,Q_T}\!\!-2\IntQ \!\!\res_f(v,\yz) e\; dxdt\\
=
\|\nabla e\|^2_{A,Q_T}
+\|\ez\|^2_{\Am,Q_T}
+\IntO e^2(x,T)dx- \|v(\cdot,0)-\xi\|^2_\Omega\\
+\int\limits_{S_{2T}}(\a\cdot n)e^2dsdt+
2\IntQ\varkappa^2 e^2dxdt-2\IntQ(e\dvg\ez+\ez\cdot\nabla e)dxdt
\end{multline}
Notice that
\ben
\label{ba6}
\IntO (e\dvg\ez+\ez\cdot\nabla e)dx=\int\limits_{S_{2T}}(\ez\cdot n)e ds
\een
and 
\ben
\label{ba7}
e^*\cdot n=(y^*\cdot n-F)=-\res_F(y^*).
\een
In view of (\ref{ba6}) and (\ref{ba7}), (\ref{ba5}) takes  the form
\begin{multline*}
\|\res_A(v,\yz)\|^2_{\Am,Q_T}\!\!-2\IntQ \!\!\res_f(v,\yz)\, e\, dxdt
+\|v(\cdot,0)-\xi\|^2_\Omega\\
=\measure^2(e,\ez)+2\|\varkappa e\|^2_{Q_T}
+\int\limits_{S_{2T}}(\a\cdot n)e^2dsdt+2\int\limits_{S_{2T}}\res_F(\yz)\,e\, dsdt,
 \end{multline*}
which implies (\ref{IDmain}).
\end{proof}
\begin{Remark} 
If $\Gamma_2=\emptyset$, $A$ is the unit matrix, $\a=0$, and $\rho=0$, then
$\varkappa=0$ and  (\ref{IDmain}) reads
\begin{multline}
\label{ba8}
\|\nabla e\|^2_{Q_T}+\|\ez\|^2_{Q_T}+\|e(\cdot,T)\|^2_\Omega\\
=
\|e(\cdot,0)\|^2_\Omega+\|\nabla v-\yz\|^2_{Q_T}-2\IntQ(f-v_t+\dvg\yz)e dxdt.
\end{multline}
This identity was earlier  derived in \cite{Re2021}. The estimate (\ref{er1})
follows from (\ref{ba8}).
\end{Remark}
Right hand side of (\ref{IDmain}) is not fully computable. It contains two directly computable norms, 
 but two other terms (integrals over $Q_T$ and $S_{2T}$)  involve unknown function $e$. 
However, these integrals can  be subject to the norms that form the
 left side of  (\ref{IDmain}). Hence the identity  generates
fully computable estimates.
We deduce them  at the end of this  section and now turn to
some special forms of (\ref{IDmain}).

%%%%%%%%%%%%%%%%%%%%%%%%%%%%%%%%%%%%%%%%%%%%
\subsection{The case $\rho>0$}
%%%%%%%%%%%%%%%%%%%%%%%%%%%
\begin{Theorem}
\label{th3}
Let the conditions (\ref{rd6})--(\ref{rd11}) be satisfied,
$\rho_\ominus>0$ (cf. (\ref{rd7})), and
$$
\rho^2_a:=\rho^2-\dvg \a\geq 0.
$$
Then, for any $v\in V_0(Q_T)$ and $\yz\in Y^*_\dvg(Q_T)$ it holds 
\begin{multline}
\label{IDad}
\measure^2(e,\ez)+\|\rho^{-1}\sigma(e^*,e)\|^2_{Q_T}+\IntQ \rho^2_a |e|^2 dxdt
+\int\limits_{S_{2T}}(\a\cdot n)e^2dsdt.\\
=
\|v(\cdot,0)-\xi\|^2_\Omega 
+\|\res_A(v,y^*)\|^2_{\Am,Q_T}
+\|\rho^{-1}\res_f(v,y^*)\|^2_{Q_T}-
2\int\limits_{S_{2T}} \res_F(y^*)\,edsdt,
\end{multline}
where 
$\sigma(e^*,e):=\dvg\ez-e_t-\a\cdot\nabla e$.

If, in addition, $\delta>0$ on $S_{2T}$, then
\begin{multline}
\label{IDadvrho}
\measure^2(e,\ez)+\|\rho^{-1}\sigma(e^*,e)\|^2_{Q_T}
+\IntQ \rho^2_a|e|^2 dxdt+\int\limits_{S_{2T}} \delta^{-1}(e^*\cdot n-\delta e)^2dsdt \\
=\|v(\cdot,0)-\xi\|^2_\Omega 
+\|\res_A(v,y^*)\|^2_{\Am,Q_T}
+\|\rho^{-1}\res_f(v,y^*)\|^2_{Q_T}+
\int\limits_{S_{2T}} \delta^{-1}\res^2_F(y^*)dsdt.
\end{multline}
\end{Theorem}
\begin{proof}
We have
\begin{multline}
\label{adv1}
\|\rho^{-1}\res_f(v,y^*)\|^2_{Q_T}=
\|\rho^{-1}(\dvg\yz+f-v_t-\a\cdot \nabla v-\rho^2  v)\|^2_{Q_T}\\
=
\|\rho^{-1}(\dvg\ez-e_t-\a\cdot\nabla e-\rho^2  e)\|^2_{Q_T}
=\|\rho^{-1}\sigma(e,e^*)-\rho  e\|^2_{Q_T}\\
=
\|\rho^{-1}\sigma(e^*,e)\|^2_{Q_T}+
\|\rho e\|^2_{Q_T}-2\IntQ\sigma(e^*,e)\, e\,dxdt.
\end{multline}

By (\ref{ba2}), (\ref{ba3}), and the relation
\be
&&\IntQ (\dvg\ez) e dxdt=-\IntQ \ez\cdot\nabla e dxdt+
\int\limits_{S_{2T}} (\ez\cdot n)\,e\,dsdt,
\ee
we find that
\begin{multline*}
2\IntQ\sigma(e^*,e) edxdt=-2\IntQ \ez\cdot\nabla e \,dxdt\\
+
2\int\limits_{S_{2T}} (\ez\cdot n)\,edsdt-\|e(x,T)\|^2_\Omega\\
+\|e(x,0)\|^2_\Omega-\int\limits_{S_{2T}}(\a\cdot n)e^2dsdt+
\IntQ (\dvg \a) \,e^2 dxdt
\end{multline*}
and modify (\ref{adv1}) as follows:
\begin{multline}
\label{adv2}
\|\rho^{-1}\res_f(v,y^*)\|^2_{Q_T}
=
\|\rho^{-1}\sigma(e^*,e)\|^2_{Q_T}
+
\|\rho e\|^2_{Q_T}
+\|e(\cdot,T)\|^2_\Omega -\|e(\cdot,0)\|^2_\Omega\\ 
+2\IntQ \ez\cdot\nabla e dxdt
-\IntQ(\dvg \a) e^2dxdt-2\int\limits_{S_{2T}} (\ez\cdot n)\,edsdt
+\int\limits_{S_{2T}}(\a\cdot n)\,e^2dsdt.
\end{multline}
Adding (\ref{adv}) and
(\ref{adv2}) we get  the identity
\begin{multline}
\label{adv4}
\measure^2(e,\ez)+\|\rho^{-1}\sigma(e^*,e)\|^2_{Q_T}
+\int\limits_{S_{2T}}(\a\cdot n)e^2dsdt+
\IntQ (\rho^2-\dvg \a) |e|^2 dxdt\\
=\|v(x,0)-\xi\|^2_\Omega+\|\rho^{-1}\res_f(v,y^*)\|^2_{Q_T}+
\|\res_A(v,\yz)\|^2_{A^{-1},Q_T}
+2\int\limits_{S_{2T}} (\ez\cdot n)\,edsdt.
\end{multline}
We recall (\ref{ba7}) and
 see that (\ref{IDad}) follows from (\ref{adv4}).

If $\delta=\a\cdot n>0$ on $S_{2T}$, then we can transform  (\ref{IDad})
 using the identity
\be
\int\limits_{S_{2T}} \left(\delta^{-1}\res^2_F(y^*)+2\res_F(y^*)e+\delta \,e^2\right)dsdt
=\int\limits_{S_{2T}} \delta^{-1}(e^*\cdot n-\delta e)^2dsdt,
\ee
Then,
\be
\int\limits_{S_{2T}}\delta\, e^2dsdt+2\int\limits_{S_{2T}}\res_F(\yz)\,e\,dsdt
=
\int\limits_{S_{2T}} \Bigl(\delta^{-1}(e^*\cdot n-\delta e)^2-
\delta^{-1}\res^2_F(y^*)\Bigr)dsdt
\ee
and
 (\ref{IDadvrho}) follows from (\ref{IDad}).
\end{proof}
\subsection{The case $\res_F(\yz)=0$}
Now we discuss corollaries of Theorems \ref{th2} and \ref{th3}
related to the case  where the term containing $\res_F(\yz)$
drops out.  This happens if
\ben
\label{ycon}
\Gamma_2=\emptyset\quad
\text{or}\quad\yz\in Y^*_F(Q_T):=\{y^*\in Y^*_\dvg(Q_T)\,\mid\,y^*\cdot n=F\;\text{on}\;S_{2T}\}.
\een
Since $\yz$ is at our disposal it is not difficult to correct this function 
so that the above condition on $S_{2T}$ is satisfied.
As a result, we obtain
identities with fully computable right hand sides.

Let the conditions of Theorem \ref{th3} be satisfied and (\ref{ycon}) holds.
 Then the last term of
(\ref{IDad}) vanishes and we arrive at the identity
\begin{multline}
\label{ID2T}
\measure^2(e,\ez)+\|\rho^{-1}\sigma(e^*,e)\|^2_{Q_T}+\IntQ \rho^2_a |e|^2 dxdt
+\int\limits_{S_{2T}}(\a\cdot n)e^2 dx\\
=
\|v(\cdot,0)-\xi\|^2_\Omega +
\|\res_A(v,y^*)\|^2_{\Am,Q_T}
+\|\rho^{-1}\res_f(v,y^*)\|^2_{Q_T}
\end{multline}
whose right hand side contains only known functions.

Let (\ref{ycon}) holds and  $\varkappa>0$. In this case, we use Theorem \ref{th2}
and rewrite 
 (\ref{IDmain})  in the form
\begin{multline}
\label{adv8}
\measure^2(e,\ez)+\|\varkappa e\|^2_{Q_T}+
\IntQ (\varkappa e+\varkappa^{-1}\res_f(v,y^*))^2dxdt\\
=
\|v(\cdot,0)-\xi\|^2_\Omega 
+\|\res_A(v,\yz)\|^2_{\Am,Q_T}
+\|\varkappa^{-1} \res_f(v,\yz)\|^2_{Q_T}. 
\end{multline}
Hence in this case,
we obtain another error identity with fully computable right hand side.

%%%%%%%%%%%%%%%%%%%%%%%%%%%%%%%%%%%%%%%%%%%%%%%%%%%
\subsection{Error estimates}
Both parts of (\ref{IDad}) and (\ref{ID2T}) contain
 integrals
with weights $\rho^{-1}$.
 Therefore, for $\rho=0$  they are not applicable.
For very small  $\rho$
these theoretically correct relations  may become
practically inefficient.
 Similar difficulties
arise in (\ref{IDadvrho}) and (\ref{adv8}) for very small $\delta$ and
$\varkappa$, respectively. In these cases, it is necessary to use
other error relations (estimates) that follow from the basic identity (\ref{IDmain}).
There are several ways to deduce
computable a posteriori estimates from (\ref{IDmain}).
Consider the simplest one. 

Let  (\ref{ycon}) hold.
Then the last integral in the right hand side
of (\ref{IDmain}) vanishes. The remaining integral
is estimated as follows:
\begin{multline*}
2\left|\int_{Q_T} \res_f(v,\yz) e\, dxdt\right|\leq
2\int\limits^T_0 \|\res_f(v,\yz)\|_\Omega\|e\|_\Omega\, dt\\
 \leq
2\frac{C_\Omega}{c_1}
\int\limits^T_0\|\res_f(v,\yz)\|_\Omega\|\nabla e\|_{A,\Omega}\,dt\\
\leq 
\int\limits^T_0
\Bigl(\frac{C^2_{A,\Omega}}{\beta^{2}(t)}\|\res_f(v,\yz)\|^2_\Omega+
\beta^2(t)\|\nabla e\|_{A,Q_T}\Bigr)dt,
\end{multline*}
where $\beta(t)$ is a positive bounded function and $C_{A,\Omega}$ is a constant in (\ref{FrA}). 
Then, for any $\beta(t)\in (0,1]$  the estimate 
\begin{multline}
\label{advinq}
\int\limits^T_0(1-\beta^2(t))\|\nabla e\|^2_{A,\Omega}dt
+\|\ez\|^2_{\Am,Q_T}+\|v(\cdot,T)-u(\cdot,T)\|^2_\Omega
+\int\limits_{S_{2T}}\delta\,e^2 dsdt\\
\leq \|v(\cdot,0)-\xi\|^2_\Omega+\|\res_A(v,\yz)\|^2_{\Am,Q_T}+
\int\limits^T_0\frac{C^2_{A,\Omega}}{\beta^{2}(t)}
 \|\res_f(v,\yz)\|^2_{\Omega} dt
\end{multline}
follows from (\ref{IDmain}).
Here (unlike (\ref{ID2T}) and (\ref{adv8})), we arrive at the estimate instead
of the identity. But it is valid for problems where $\rho$ and $\varkappa$ are
 small or even equal to zero.
 
A lower bound is derived by the analogous procedure. It reads
\begin{multline}
\label{advmin}
\int\limits^T_0(1+\beta^2(t))\|\nabla e\|^2_{A,\Omega}dt
+\|\ez\|^2_{\Am,Q_T}+\|v(\cdot,T)-u(\cdot,T)\|^2_\Omega
+\int\limits_{S_{2T}}\delta\,e^2 dsdt\\
\geq \|v(\cdot,0)-\xi\|^2_\Omega+\|\res_A(v,\yz)\|^2_{\Am,Q_T}-
\int\limits^T_0\frac{C^2_{A,\Omega}}{\beta^{2}(t)}
 \|\res_f(v,\yz)\|^2_{\Omega} dt.
\end{multline}
Certainly this bound is meaningful only if the right hand side of (\ref{advmin}) is positive.

%%%%%%%%%%%%%%%%%%%%%%%%%%%%%%%%%%%%%%%%%%%%%%%%%%%%%%%%%%%%%%%%%%
\section{Applications to analysis of inverse problems}
The estimates derived above 
can be applied not only to the  problem 
(\ref{rd1})--(\ref{rd5}). In this section, we consider
several problems  based on the equation (\ref{rd1}) with different 
conditions imposed on solution and problem data, which
originate from various real life
models  in experimental sciences and engineering.

As in (\ref{er2})--(\ref{er3}), 
  $v_\eps\in V_0(Q_T)$ (cf. (\ref{Vnul})) and $\yz_\eps\in Y^*_\dvg(Q_T)$ (cf. (\ref{Ydiv})).
 These functions may
not exactly satisfy the relations (\ref{rd1}) and (\ref{rd2}), but they satisfy the boundary
conditions of an ill--posed problem. In particular, $v_\eps$ and $\yz_\eps$ can be
computed by 
  some
numerical procedure developed for the corresponding class of inverse problems. 
Theorems \ref{th2} and \ref{th3}
generate functionals (e.g., (\ref{Ubet})) that control the distance between $(v_\eps,\yz_\eps)$
and a pair  $(u_\eps,\pz_\eps)$ formed by exact solutions
of a direct boundary value problem associated with the same differential equation
with the boundary conditions defined by  $(v_\eps,\yz_\eps)$.
Assume that this functional computed for $(v_\eps,\yz_\eps)$ is smaller than $\eps^2$. In this case,   
we can guarantee that 
$(u_\eps,\pz_\eps)$ and $(v_\eps,\yz_\eps)$ are close pairs of functions. Moreover, $(u_\eps,\pz_\eps)$
inherit some of the boundary conditions from  $(v_\eps,\yz_\eps)$. 
In other words, they exactly satisfy one part
of the  required conditions  while the others hold approximately with the accuracy $\eps$.
We note that there is no single, definitive choice here: the conditions to be satisfied
exactly and those only approximately may
be defined in various ways (see examples below).
As a result, we find a well-posed problem, the solution of which 
can be considered as an
$\eps$- approximation for the ill-posed one.

%%%%%%%%%%%%%%%%%%%%%%%%%%%%%%%%%%%%%
\subsection{Backward  problem}
 First, we consider the problem: find $ u$ and $ \pz$ such that
\ben
\label{bw1}
&& u_t-\dvg  \pz+\a\cdot\nabla  u+\rho^2   u=f\;{\rm in}\;Q_T,\\
\label{bw2}
&& \pz=A\nabla  u,\\
\label{bw3}
&& u(x,t)=0\qquad{\rm on}\;S_{1T},\qquad  \pz\cdot n=F\;\text{on}\;S_{2T},\\
\label{bw4}
&& u(x,T)=\chi(x)\in L_2(\Omega).
\een
Problem (\ref{sim1})--(\ref{sim3}) is a simple particular case of this
backward parabolic problem.
%COMMENT of DeepSeek: The direct problem is posed with $\xi=v_\eps(x,0)$, but $v_\eps\in V_0(\Q)$
%only guarantees $v_\eps(\cdot,0)\in L_2(\Om)$ by trace. The assumption
%\eqref{rd11} requires $\xi\in\mathring W^1_2(\Om)$. Either strengthen the
%assumption on $v_\eps$ (e.g.\ $v_\eps\in\mathring W^{1,1}_2(\Q)$ with
%$v_\eps(\cdot,0)\in\mathring W^1_2(\Om)$) or invoke a density/limiting
%argument.
% *** ANSWER: This is an ill posed problem. It is never claimed that it has a solution
% with so defined conditions.

For $\beta\in (0,1]$ and $\gamma>0$, we define the functional
\ben
\label{Ubet}
\Upsilon_{\beta,\gamma}(v,\yz):=\|\res_A(v,\yz)\|^2_{\Am,Q_T}+
\int\limits^T_0 \frac{C^2_{A,\Omega}}{\beta^{2}(t)}
 \|\res_f(v,\yz)\|^2_{\Omega} dt+\frac{1}{\gamma}\|v(\cdot,T)-\chi\|^2_\Omega.
\een

\begin{Theorem}
\label{th4}
Let $v_\eps\in V_0(Q_T)$ and $\yz_\eps\in Y^*_\dvg(Q_T)$ satisfy the condition
$
\yz_\eps\cdot n=F\;\text{on}\;\Gamma_{2T},
$
and
 \ben
 \label{bwUbet}
 \Upsilon_{\beta,\gamma}(v_\eps,y^*_\eps)\leq \eps^2.
 \een
Then, problem (\ref{rd1})--(\ref{rd6}) with the initial condition $\xi=v_\eps(x,0)$
has the solution $u_\eps$ such that
\ben
\label{eps1}
&&\|u_\eps(x,T)-\chi(x)\|_\Omega\leq\eps(1+\gamma)^{1/2}.
\een
The solution $u_\eps$ and its flux $\pz_\eps$ belong to $\eps$-neighbourhoods of 
$v_\eps$ and $\yz_\eps$, respectively, i.e., 
\ben
\label{eps2}
&&\|\pz_\eps-\yz_\eps\|_{Q_T}\leq c_2\eps,
\een

and
\ben
\label{eps3}
\int\limits^T_0(1-\beta^2)\|\nabla(u_\eps-v_\eps)\|^2_{A,\Omega}dt+
\int\limits_{S_{2T}}\delta |u_\eps-v_\eps|^2 dsdt
\leq \eps^2.
\een
\end{Theorem}
\begin{proof}
We apply (\ref{advinq}) to
the correct initial boundary value problem (\ref{rd1})--(\ref{rd6})
with the initial condition $\xi=v_\eps(x,0)$ and consider $v=v_\eps$ and $\yz=\yz_\eps$
as approximations of the corresponding exact functions $u_\eps$ and $\pz_\eps:=A\nabla u_\eps$.
We have
\begin{multline*}
\int\limits^T_0(1-\beta^2(t))\|\nabla (v_\eps-u_\eps)\|^2_{A,\Omega}dt
+\|\yz_\eps-\pz_\eps\|^2_{\Am,Q_T}\\
+\|v_\eps(\cdot,T)-u_\eps(\cdot,T)\|^2_\Omega+
\frac{1}{\gamma}\|v_\eps(\cdot,T)-\chi\|^2_\Omega\\
+\int\limits_{S_{2T}}\delta\,|v_\eps-u_\eps|^2 dsdt
\leq \Upsilon_{\beta,\gamma}(v_\eps,\yz_\eps)\leq \eps^2.
\end{multline*}
Since 
\be
\|v_\eps(\cdot,T)-u_\eps(\cdot,T)\|^2_\Omega\geq 
\frac{1}{1+\gamma}\|u_\eps(\cdot,T)-\chi\|^2_\Omega-\frac{1}{\gamma}
\|v_\eps(\cdot,T)-\chi\|^2_\Omega 
\ee
for any $\gamma>0$, we find that
\begin{multline}
\label{advbet}
\int\limits^T_0(1-\beta^2(t))\|\nabla (v_\eps-u_\eps)\|^2_{A,\Omega}dt
+\|\yz_\eps-\pz_\eps\|^2_{\Am,Q_T}\\
+\!\int\limits_{S_{2T}}\delta |v_\eps-u_\eps|^2 dsdt+\frac{1}{1+\gamma}\|u_\eps(\cdot,T)-\chi\|^2_\Omega
\leq \eps^2,
\end{multline}
%%% *** NO, your comment on the term $\|\yz_\eps-\pz_\eps\|^2_{\Am,\Omega}$
%%%       seems to me superfluous. The above norn is a function of $t$ and the factor
%%%.  $(1-\beta(t))$ is acting on the first term only. 
wherefrom  the estimates (\ref{eps1}), (\ref{eps2}), and (\ref{eps3}) follow.
\end{proof}

\begin{Remark}
Problem (\ref{sim1})--(\ref{sim3}) is a very special case of (\ref{bw1})--(\ref{bw4}).
It is easy to see that (\ref{er2}) follows from (\ref{advbet}) if $\a=0$, $\rho=0$, $S_{2T}=\emptyset$,
and $A$ is the unit matrix.  If the condition $v_\eps(x,T)=\chi(x)$ is exactly satisfied, then the last term in (\ref{Ubet})
vanishes and we can set $\gamma=0$. Then (\ref{Ubet}) is reduced to
$\Upsilon_\beta(u_\eps,\yz_\eps)$ and (\ref{bwUbet}) to (\ref{Iopt}).
\end{Remark}

If $\rho>0$ or $\varkappa>0$, then modified versions of Theorem
\ref{th4}  follow from   (\ref{ID2T})
and (\ref{adv8}).
Assume that $\rho>0$, (\ref{ycon}) holds,  and 
\be
\Upsilon_\rho(v_\eps,y^*_\eps):=\; 
\|\res_A(v_\eps,y^*_\eps)\|^2_{\Am,Q_T}
+\|\rho^{-1}\res_f(v_\eps,y^*_\eps)\|^2_{Q_T}\leq \eps^2
\ee 
for $v_\eps\in V_0(Q_T)$ and $\yz_\eps\in Y^*_\dvg(Q_T)$ and, in addition, $v_\eps(x,T)=\chi(x)$.

By (\ref{ID2T}) and quite analogous arguments, we prove that direct
problem (\ref{rd1})--(\ref{rd6}) with $\xi(x)=v_\eps(x,0)$ has the solution
$(u_\eps,\pz_\eps)$ and
\ben
\label{advrho}
\measure^2(e,\ez)+\|\rho^{-1}\sigma(e^*,e)\|^2_{Q_T}
+\IntQ \rho^2_a |e|^2 dxdt+\int\limits_{S_{2T}}\delta \,e^2dsdt\;\leq\;\eps^2,
\een
where $e=v_\eps-u_\eps$ and $\ez=\yz_\eps-\pz_\eps$. Hence 
\ben
\label{adveps}
\|\nabla (u_\eps-v_\eps)\|_{Q_T}\leq \eps,\quad \|\pz_\eps-\yz_\eps\|_{Q_T}\leq\eps,\quad
\text{and}\;\|u_\eps(x,T)-\chi(x)\|_\Omega\leq\eps.
\een
We conclude that $\eps$--neighbourhood of $(v_\eps,y^*_\eps)$ contains
exact solution of a correct initial boundary--value problem with the initial
condition $\xi=v_\eps(x,0)$.
 \vskip3pt
Similar reasoning can be used if $\varkappa>0$. 
Let 
$$v_\eps\in V_0(Q_T),\;\;  v_\eps(x,T)=\chi(x),\;\; \text{and}\; \;\yz_\eps\in Y^*_\dvg(Q_T)
$$ 
be such that 
\be
\Upsilon_\varkappa(v_\eps,y^*_\eps):=\,\|\res_A(v_\eps,y^*_\eps)\|^2_{\Am,Q_T}
+\|\varkappa^{-1}\res_f(v_\eps,y^*_\eps)\|^2_{Q_T}\,\leq \,\eps^2.
\ee
Then, by (\ref{adv8}) we conclude that
\ben
\label{advkap}
\measure^2(e,\ez)+\|\varkappa e\|^2_{Q_T}\;\leq \eps^2
\een
and, therefore, the estimates (\ref{adveps}) hold.

\begin{Remark}
Estimate (\ref{advbet}) shows that the sum of the four nonnegative 
terms is less than $\eps^2$. The
quantity $\|\chi(x)~-~u_\eps(x,T)\|^2_\Omega$ is only one of these terms. Therefore,  we can
expect that the norm $\|\chi(x)~-~u_\eps(x,T)\|_\Omega$ is essentially smaller than $\eps$.
Analogous conclusions follow from (\ref{advrho}) and (\ref{advkap}).
In other words, it is very likely that $u_\eps$   approximates  the desired 
condition at $t=T$ much better. Since $u_\eps(x,0)$ is known,
we may try to find a sharp approximation of $u_\eps$.
Stable and efficient numerical 
methods for solving
direct initial boundary value problems are well developed
(e.g. see \cite{Thomee}). Using them in combination with respective
 a posteriori error estimation methods (discussed in Sect. \ref{ApoIdn}) suggests
 a way to get efficient and controllable approximations for the backward
 problem (\ref{bw1})--(\ref{bw4}) with the confirmed accuracy
 wich is much better  than
 it follows from (\ref{adveps}).
\end{Remark}
%%%%%%%%%%%%%%%%%%%%%%%%%%%%%%%%%%%%%%%%%%%%%%%%%%%%%%%%%%%%%%%%%%
\subsection{Problem with nonlocal conditions}
Consider the problem 
\ben
\label{nl1}
&&u_t-\dvg \pz+\a\cdot\nabla u+\rho^2  u=f\;{\rm in}\;Q_T,\\
\label{nl2}
&&\pz=A\nabla u,\\
\label{nl3}
&&u(x,t)=0\qquad{\rm on}\;S_{T},\\
\label{nl4}
&&u(x,0)+\omega\,u(x,T)=0,
\een
where $\omega\in [-1,0)$ is a given function.  
Existence and uniqueness of
such type problems have been studied (e.g., see \cite{Kerefov,Chab}).

Assume that we have the function $v_\eps\in V_0(Q_T)$ satisfing the condition
\ben
\label{nl5}
&&v_\eps(x,0)+\omega\, v_\eps(x,T)=0
\een
and $\yz_\eps\in Y^*_\dvg(Q_T)$ such that for $\beta(t)\in (0,1]$
\ben
\label{nl6}
\Upsilon_{\beta}(v_\eps,\yz_\eps):=\|\res_A(v_\eps,\yz_\eps)\|^2_{\Am,Q_T}+
\int\limits^T_0\frac{C^2_{A,\Omega}}{\beta^{2}(t)}
 \|\res_f(v_\eps,\yz_\eps)\|^2_{A,\Omega} dt\leq \eps^2.
\een

Let $u_\eps(x,t)$ and $\yz_\eps$ solve the
  initial boundary value problem
  \ben
\label{nl1eps}
&&(u_\eps)_t-\dvg \pz_\eps+\a\cdot\nabla u_\eps+\rho^2  u_\eps=f\;{\rm in}\;Q_T,\\
\label{nl2eps}
&&\pz_\eps=A\nabla u_\eps,\\
\label{nl3eps}
&&u_\eps(x,t)=0\qquad{\rm on}\;S_{T},\\
\label{nl4eps}
&&u_\eps(x,0)=v_\eps(x,0).
\een 
This is a well posed problem. 
We consider $v_\eps$ and $y_\eps$ as approximations
of the well posed problem (\ref{nl1eps})--(\ref{nl4eps})
and use (\ref{advinq}) with $\xi(x)=v_\eps(x,0)$ (note that in our case $S_{2T}=\emptyset$). 
By (\ref{advinq}) and (\ref{nl6}),
we find that
\begin{multline}
\label{nl5eps}
\int\limits^T_0(1-\!\beta^2(t))\|\nabla(v_\eps-u_\eps)\|^2_{A,\Omega}dt\\
+\|\yz_\eps-\pz_\eps\|^2_{\Am,Q_T}+\|v_\eps(x,T)-u_\eps(x,T)\|^2\,\leq \,\eps^2.
\end{multline}
Since
\be
\|v_\eps(x,T)-u_\eps(x,T)\|_\Omega=\|\omega^{-1}v_\eps(x,0)+u_\eps(x,T)\|_\Omega
=|\omega|^{-1}\| u_\eps(x,0)+\omega u_\eps(x,T)\|_\Omega,
\ee
from (\ref{nl5eps}) we deduce the estimate 
\begin{multline*}
\int\limits^T_0(1-\beta^2(t))\|\nabla (v_\eps-u_\eps)\|^2_{A,\Omega}dt
+\|\yz_\eps-\pz_\eps\|^2_{\Am,Q_T}\\
+|\omega|^{-2}\|u_\eps(x,0)+\omega u_\eps(x,T)\|^2\leq \eps^2.
\end{multline*}
Thus, we arrive at the following result.
\begin{Theorem}
\label{th5}
Let $v_\eps\in V_0(Q_T)$ and $\yz_\eps\in Y^*_\dvg(Q_T)$ satisfy the conditions
(\ref{nl5}) and (\ref{nl6}).
Then, solution of (\ref{nl1eps})--(\ref{nl4eps}) meets the estimate
\ben
\label{nl6eps}
&&\|u_\eps(x,0)+\omega u_\eps(x,T)\|_\Omega\leq\;|\omega|\eps.
\een
Moreover,
 $u_\eps$ and  $\pz_\eps$ belong to $\eps$--neighbourhoods of 
$v_\eps$ and $\yz_\eps$, respectively, i.e., 
\ben
\label{nl7eps}
&&\|\pz_\eps-\yz_\eps\|_{Q_T}\leq \eps,
\een
and
\ben
\label{nl8eps}
\int\limits^T_0(1-\beta^2)\|\nabla(u_\eps-v_\eps)\|^2_{A,\Omega}\,dt\;
\leq\; \eps^2.
\een
\end{Theorem}
Theorem \ref{th5} has a clear meaning. In view of
(\ref{nl7eps}) and (\ref{nl8eps}), we guarantee that $\eps$-neighbourhoods of 
$v_\eps$ and $\yz_\eps$ contain
solutions of a direct problem with the initial condition $\xi=v_\eps(x,0)$
such that relations (\ref{nl1}) -- (\ref{nl3}) are exactly satisfied
 and (\ref{nl6eps}) shows that
(\ref{nl4}) is satisfied approximately with the accuracy $\omega\eps$.

%%%%%%%%%%%%%%%%%%%%%%%%%%
\subsection{Inverse Cauchy problem}
Inverse Cauchy problem imposes
an over-specified boundary condition on one part of the boundary
while the other part  remains free from any conditions.
Usually, these problems are considered for
elliptic equations but conditions of this type also
generate a class of inverse parabolic problems (e.g., see Chapt. 5 of
\cite{LattLion}). Modifications of these problems are studied
in connection with various applied problems. At this point we address the
reader to 
\cite{Engl,KabLav,SamVab,MiM} and numerous publications 
 cited therein.
 
As an example, we take the problem
\ben
\label{ca1}
&u_t-\dvg \pz+\a\cdot\nabla u+\rho^2  u=f\;&{\rm in}\;Q_T,\\
\label{ca2}
&\pz=A\nabla u,\\
\label{ca3}
&u(x,t)=0\qquad&{\rm on}\;S_{2T},\\
\label{ca4}
& \pz\cdot n= F\qquad&{\rm on}\;S_{2T},\\
\label{ca5}
&u(x,0)=\xi(x).
\een
Here two boundary conditions are stated on one part of $S_T$
(let it be $S_{2T}$) and the Dirichlet boundary conditions on the remaining part ($S_{1T})$ are to be defined
\begin{Theorem}
\label{th6}
Let $v_\eps\in V$
and $\yz_\eps\in Y^*_\dvg(Q_T)$ satisfy the conditions
\be 
 v_\eps(x,t)=0\; \text{\rm on}\; \Gamma_{2T}, \qquad v_\eps(x,0)=\xi(x),\qquad
\yz_\eps\cdot n~=~ F\;\text{\rm on}\;S_{2T},
\ee 
and
 \ben
 \label{caUbet}
 \Upsilon_\beta(v_\eps,y^*_\eps)\leq \eps^2\quad\text{for}\quad \beta(t)\in(0,1].
 \een

Then,
\ben
\label{ep1}
&&\int\limits^T_0\!(1\!-\!\beta^2)\|\nabla(u_\eps-v_\eps)\|^2_{A,\Omega}\leq \eps^2,\;\;\;\;
\|\pz_\eps-\yz_\eps\|_{Q_T}\leq  c_2\eps,\\
\label{ep2}
&&\|u_\eps(x,T)-v_\eps(x,T)\|_\Omega\leq\eps,\\
\label{ep3}
&&\|u_\eps\|_{\Gamma_{2T}}\leq C_{\Omega,\Gamma_2}
\left(1-\beta^2_\oplus \right)^{-1/2}\eps,
\een
where $\beta_\oplus=\sup\limits_{t\in [0,T]}\beta(t)$,
$C_{\Omega,\Gamma_2}$ is a constant in the  trace
inequality associated with $\Gamma_2$, and
$u_\eps$ and $\pz_\eps$ solve the initial boundary value problem
  \ben
\label{ch1eps}
&(u_\eps)_t-\dvg \pz_\eps+\a\cdot\nabla u_\eps+\rho^2  u_\eps=f\;&{\rm in}\;Q_T,\\
\label{ch2eps}
&\pz_\eps=A\nabla u_\eps,\\
\label{ch3eps}
&u_\eps(x,t)=v_\eps(x,t)\;&{\rm on}\;S_{1T},\\
\label{ch4eps}
&\pz_\eps\cdot n=  F\;\; &\text{\rm on}\; S_{2T},\\
\label{ch5eps}
&u_\eps(x,0)=\xi(x).
\een 
\end{Theorem}
\begin{proof}
We apply
(\ref{advinq}) to the problem (\ref{ch1eps})--(\ref{ch5eps}).
 In view of (\ref{caUbet}), we have
\be
\int\limits^T_0(1\!-\beta^2(t))\|\nabla (u_\eps-v_\eps)\|^2_{A,\Omega}dt
+\|\pz_\eps-\yz_\eps\|^2_{\Am,Q_T}
+\|v_\eps(\cdot,T)-\!u_\eps(x,T)\|^2_\Omega
\leq \eps^2,
\ee
wherefrom (\ref{ep1}) and (\ref{ep2}) follow. Since
\ben
\label{ch6eps}
\|u_\eps\|_{\Gamma_{2T}}\leq \|u_\eps-v_\eps\|_{\Gamma_{2T}}+\|v_\eps\|_{\Gamma_{2T}}\leq
C_{\Omega,\Gamma_2}\|\nabla(u_\eps-v_\eps)\|_\Omega,
\een
 (\ref{ep3}) follows from the first inequality in (\ref{ep1})
 and (\ref{ch6eps}). 
\end{proof}
Theorem  \ref{th6} guarantees that $\eps$-neighbourhood of known functions $v_\eps$ and $\yz_\eps$
contains exact solution of the well-posed problem (\ref{ch1eps})--(\ref{ch5eps}), which
satisfies the equation (\ref{ca1}), relation (\ref{ca2}) 
and conditions (\ref{ca4}) and (\ref{ca5}). The condition (\ref{ca3}) is satisfied
with $\eps$-accuracy.
\begin{Remark}
In \cite{RepinCMAM2026},   ill posed Cauchy problems for elliptic boundary value
problems are analyzed in the framework of the above discussed 
approach based on functional a posteriori estimates.
\end{Remark}

%%%%%%%%%%%%%%%%%%%%%%%%%%%%%%%%%%%%%%%%%%%%%%%%%%%%%%%%%%%%%%%%%%%%%%%%%%%%%%%%%%%
\subsection{Identification problems}
This class of problems is motivated by a situation often encountered in applications where
model parameters must be determined from a set of available
data (e.g., see \cite{BanKun,Kab}). In the context of the  initial boundary value problem under consideration, 
such unknown parameters
may be  $f$, $A$, $\a$, and $\rho$.

For example, if $f$ is sought and the other parameters are given, 
then we are dealing 
with a "source reconstruction" problem. 
The methodology discussed above is entirely 
applicable to this class of problems.
Let us briefly consider this issue with the paradigm of the problem:
find $f(x,t)$ such that the relations (\ref{nl1})--(\ref{nl3}) hold and,
in addition,
\ben
\label{par1}
&&u(x,0)=\xi(x)\qquad\text{\rm and}\qquad u(x,T)=\chi(x).
\een
Here $\xi$ and $\chi$ are certain given functions and $f$
is a function in a set $U_f$ of "admissible"
source terms. Depending on how the set $U_f$ is defined, this problem may or may 
not have a solution.

Assume that we know $f_\eps\in U_f$ (e.g., $U_f=L_1(0,T;L_2(\Omega))$), $v_\eps\in V_0(Q_T)$, and $\yz_\eps\in Y^*_\dvg(Q_T)$
such that $v_\eps(x,0)=\xi$, $v_\eps(x,T)=\chi$,
and
\ben
\label{paridn}
\Upsilon_\beta(v_\eps,\yz_\eps)\leq \eps^2.
\een
Direct problem with $f=f_\eps$  and  initial condition $u(x,0)=\xi$ is well
posed and has a solution $u_\eps$. 
We apply (\ref{advinq}) (note that $\Gamma_2=\emptyset$),
use (\ref{paridn}), and obtain
\ben
\label{parbet}
\int\limits^T_0\!\Bigl((1-\!\beta^2(t))\|\nabla (v_\eps-u_\eps)\|^2_{A,\Omega}+
\|\yz_\eps-\pz_\eps\|^2_{\Am,\Omega}\Bigr)dt+\|\chi(x)\!-u_\eps(x,T)\|^2_\Omega
\leq \eps^2.
\een
Hence we can define neighbourhoods of $v_\eps$
and $\yz_\eps$ (which are proportional to $\eps$) that contain $u_\eps(x,t)$ and
 $\pz_\eps=A\nabla u_\eps(x,t)$, respectively. The function 
 $u_\eps(x,t)$ exactly satisfies the relations
(\ref{nl1}) -- (\ref{nl3}) and the initial condition $u_\eps(x,0)=\xi(x)$.
The condition at $t=T$ is satisfied approximately:
 $\|\chi(x)\!-u_\eps(x,T)\|_\Omega\leq\eps$.

This scheme is also suitable for more complex identification 
problems where a larger number of conditions are imposed 
on the solution and a larger number of parameters are 
to be determined based on them. For example, the conditions 
on $A$ can be weakened in the sense that it is not a given matrix, but an element 
of some set of admissible matrices $U_A$ (that satisfy (\ref{rd10})).
Analogously, it may be known that $\a\in U_\a$ and $\rho\in U_\rho$.
%Additional conditions 
%$$u(x,t_k)=\chi_k(x)\qquad k=1,2,...,N
%$$
%may be given at certain   intermediate
%points $t_k$ inside the interval $(0,T)$.

Assume that we know the functions
$$
(f_\eps,A_\eps,\a_\eps,\rho_\eps)\in U:=U_f\times U_A\times U_\a\times U_\rho
$$
that define a particular initial boundary value problem
 and the functions $v_\eps\in V_0(Q_T)$ and $\yz_\eps\in Y^*_\dvg(Q_T)$
 (viewed as approximations for this problem)
such that
\ben
\label{Uident}
\Upsilon_\beta(v_\eps,\yz_\eps;f_\eps,A_\eps,\a_\eps,\rho_\eps)\leq \eps^2.
\een
The estimate (\ref{Uident}) means that  $\eps-$neighbourhoods of $v_\eps$
and $\yz_\eps$ contain $u_\eps$ and $\pz_\eps$ that present the solution
of the direct evolutionary problem
with the initial condition $\xi=v_\eps(x,0)$ and data
$f_\eps$, $A_\eps$,
$\a_\eps$, and $\rho_\eps$.
%%%%%%%%%%%%%%%%%%%%%%%%%%%%%%%%%%%%%%%%%%%%%  
\subsection{Applications to numerical methods}
\label{num}
The question of how to find functions $u_\eps$ and $\yz_\eps$
 is a separate large topic which is not discussed 
in this article. At this point, we refer to numerous publications 
(e.g., see \cite{Engl,Kab,KabLav,SamVab,LattLion,TikhAr})
devoted to numerical analysis of ill-posed problems.
It is worth noting that results of Section 4 
also suggest a practical way to solve the
 problems discussed in this section. 

In the above, the
 functionals $\Upsilon_\beta(v,\yz)$, $\Upsilon_\rho(v,\yz)$
and a more complicated functional $\Upsilon_{\beta,\gamma}$ (defined by (\ref{Ubet}))
 are used to
verify the accuracy of approximate solutions.
However, these functionals
can also be considered as the basis for numerical methods 
aimed to find an approximate solution.
For example, consider  the problem (\ref{bw1})--(\ref{bw4}). 
 It is easy to see that
finding suitable $v_\eps$ and $\yz_\eps$ is reduced to the quadratic
minimization problem
\ben
\label{minpro}
\inf\limits_{v\in V_0(Q_T),\atop y^*\in Y^*_\dvg(Q_T)} \Upsilon_\beta(v,y^*).
\een
In (\ref{minpro}),
a minimizing sequence $(v_k,\yz_k)$, $k=1,2,...$ can be computed by known numerical
methods.  Theorem \ref{th4} claims that the value of
$\Upsilon_\beta(v_k,y^*_k)$ shows guaranteed
 accuracy of the approximation achieved on the step $k$. 
 
 For backward parabolic problems, we may
 await that infimum in (\ref{minpro}) is equal to zero (cf. (\ref{bwmod3})).
 In this case, it is theoretically possible to find $v_k,y^*_k$ such that $\Upsilon_\beta(v_k,y^*_k)$
 is arbitrarily small (what, however, may be not an easy task). The situation 
 is drastically different
 for parameter identification problems, where the distance to the set of
 exact solutions of "admissible" direct problems may be finite so that
  $\eps$ in (\ref{paridn}) or (\ref{Uident}) cannot be made smaller than
 some positive value. Here the  notion of "solution" should be defined in the Tikhonov's
 sense and the value of  $\Upsilon_\beta(v_k,y^*_k)$ gives an upper bound of the
 distance to the set formed by solutions of "close" well-posed problems
 that
 contains a "quasi-solution".
%%%%%%%%%%%%%%%%%%%%%%%%%%%%%%%%%%%%%%%%%%%%%%%


\begin{thebibliography}{363}

\bibitem{BanKun}
H. T. Bank, K. Kunisch.  Estimation Techniques for Distributed Parameter System. 
 Birkh\"auser, Boston, 1989.

\bibitem{Chab}
J. Chabrowski. On non-local problems for parabolic equations. Nagoya Math. J.
Vol. 93 (1984), 109-131.

\bibitem{Engl}
H. W. Engl, Inverse problems and their regularization. In: Computational Mathematics Driven
by Industrial Problems. Springer, Berlin–New York, 1999, 126–150.

\bibitem{Kab}
S. I. Kabanikhin
Definitions and examples
of inverse and ill-posed problems.
J. Inv. Ill-Posed Problems 16 (2008), 317–357.

\bibitem{KabLav}
S. I. Kabanikhin and M. M. Lavrentiev (Eds), Ill-Posed and Non-Classical Problems of Mathematical
Physics and Analysis. VSP, Utrecht, 2003.

\bibitem{LattLion}
R. Lattes and J.-L. Lions, M\'ethode de Quasi Reversibilit\'e et Applications.
 Dunod, Paris, 1967.

\bibitem{SamVab}
A. A. Samarskii and P. N. Vabishchevich, Numerical Methods for Solving Inverse Problems of
Mathematical Physics. VSP, Utrecht, 2007.

\bibitem{TikhAr}
A. N. Tikhonov and V. Ya. Arsenin, Solutions of Ill-Posed Problems. Wiley, New York, 1977.


\bibitem{Kerefov}
A. A. Kerefov. 
Nonlocal boundary value problems for parabolic equations.
Differentsial'nye Uravneniya, 1979, Volume 15, Number 1, Pages 74–78 (in Russian).

%%%%%%%%%%%%%%%%%%%%%%%%%%%%%%%%%%%%%%%%%%%%%%%%%%%%
\bibitem{LaSoUr}
O. A. Ladyzhenskaya, V. A. Solonnikov and N. N. Uraltseva.
{\it Linear and Quasilinear Equations of Parabolic Type},
Nauka, Moscow, 1967.

%%%%%%%%%%%%%%%%%%%%%%%%%%%%%%%%%%%%%%%%%%%%%%%%%%%%
\bibitem{La}
O.~A.~Ladyzhenskaya,
{\it The boundary value problems of mathematical physics},
Springer, New York, 1985.

\bibitem{MiM}
A. S. Mikhaylov, V. S. Mikhaylov, Discrete dynamical
systems, inverse problems and related topics. Journal of Inverse
and Ill-posed Problems 10.1515/jiip-2024-0062, 23 (2), 217--253, 2025.
%%%%%%%%%%%%%%%%%%%%%%%%%%%%%%%%%%%%%%%%%%%%%%%%%%%%%
\bibitem{ReLincei}
S. Repin. Estimates of deviations from exact solutions initial-boundary value problem for the heat equation. Atti Accad. Naz. Lincei Cl. Sci. Fis. Mat. Natur. Rend. Lincei (9) Mat. Appl., 13(2):121--133, 2002.

\bibitem{ReGruyter}
S. Repin. A posteriori estimates for partial differential equations, 
Walter de Gruyter GmbH \& Co. KG, Berlin, 2008.

\bibitem{Re2021}
S. Repin.
Error identities for parabolic initial
boundary value problems. 
Zapiski Nauchn. Semin. Steklov Math. Inst. (POMI).
%Записки научных семинаров ПОМИ, 
508  (2021), 
147--172; J. Math. Sci., 283(2024),797–815.

\bibitem{Re2022}
S. Repin. A posteriori error identities for parabolic convection-diffusion problems.
Zapiski Nauchn. Semin. Steklov Math. Inst. (POMI), 519(2022), 205--228;
J. Math. Sci., 297 (2026), 606–622. 

\bibitem{RepinCMAM2026}
S. Repin. A posteriori analysis of inverse problems for elliptic differential equations.
Computational Methods in Applied Mathematics, 4(26), 2026.

 \bibitem{Show}
R. E. Showalter. The Final Value Problem 
 for Evolution Equations.
  Journal of Mathematical Analysis and Applications 47 (1974), 563--572.
 
\bibitem{Thomee}
V. Thom\'ee. Galerkin finite element methods for parabolic problems.
Springer Series in Computational Mathematics, 25, Springer, Berlin, 2006.
\end{thebibliography}
\end{document}